\documentclass[smallextended,numbook,runningheads]{svjour3}
\usepackage{amsmath}
\smartqed  

\usepackage{amsfonts,amssymb}
\usepackage{epsfig}
\usepackage{booktabs}
\usepackage{xcolor}

\journalname{}
\date{ \phantom{b} \vspace{45mm}\phantom{e}}
\def\makeheadbox{\hspace{-4mm}Version of 2 November 2022}

\def\real{{\mathbb R}}

\def\d{{\mathrm d}}
\def\e{{\mathrm e}}

\def\iu{\mathrm{i}}
\def\eps{\varepsilon}

\def\forall{\qquad\hbox{ for all }\ }

\newdimen\GGGlength
\newdimen\GGGheight
\newbox\GGGbox
\def\GGGput[#1,#2](#3,#4)#5{%
  \setbox\GGGbox\vbox{\hbox{#5}\kern0pt}%
  \GGGlength\wd\GGGbox%
  \divide\GGGlength by100 \multiply\GGGlength by#1%
  \GGGheight\ht\GGGbox%
  \divide\GGGheight by100 \multiply\GGGheight by#2%
  \put(#3,#4){\kern-\GGGlength\raise-\GGGheight\box\GGGbox}}

\begin{document}

\title{Drift approximation by the modified Boris algorithm of charged-particle dynamics in toroidal geometry}

\titlerunning{Drift approximation of charged-particle dynamics in toroidal geometry}

\author{
Yanyan~Shi$^{1}$}
\authorrunning{
Y.\ Shi}

\institute{
$^1$~Mathematisches Institut, Univ.\ T\"ubingen, D-72076 T\"ubingen, Germany.\\
\phantom{$^1$~}\email{shi@na.uni-tuebingen.de}\\
}

\date{ }

\maketitle

 \begin{abstract} 
 In this paper, we  study the charged-particle dynamics under strong magnetic field in a toroidal axi-symmetric geometry. 
 Using modulated Fourier expansions of the exact and numerical solutions, the long-term drift motion of the exact solution in toroidal geometry  is derived and  the error analysis of the large-stepsize modified Boris algorithm over long time scales is provided.
Numerical experiments illustrate the theoretical results.

 \bigskip

\end{abstract}

\section{Introduction}
The time integration of the equations of motion of charged particles is a crucial step in particle methods of plasma physics~\cite{birdsall05ppv}. In the strong magnetic field regime, the charged particles exhibit very fast rotations of small radius around a guiding centre. This often brings stringent restriction on the time stepsize for numerical integrators. There are many works aiming at designing large-stepsize integrators with good accuracy for the charged-particle dynamics such as~\cite{vu95ans,filbet2017asymptotically,chartier2020uniformly,ricketson20aec,xiao21smc,hairer2022large}. Among them, a Boris-type integrator with appropriate modifications shows striking numerical results~\cite{xiao21smc} and the rigorous analysis is provided in~\cite{lubich2022large}. It is proved that the position and the parallel velocity are approximated with $O(h^2)$ accuracy for the modified Boris algorithm with large step sizes $h^2\sim \eps$ for fixed $T=O(1)$, where $\eps\ll 1$  is a small parameter whose inverse corresponds to the strength of the magnetic field.

In this paper, we are interested in analyzing the long time behavior (over O($\eps^{-1}$)) of the modified Boris algorithm in a toroidal axi-symmetric geometry, with a magnetic field everywhere  toroidal and an electric field everywhere orthogonal to the magnetic field. This geometry has already been proposed in~\cite{filbet2020asymptotics} and a first order description of the slow dynamics for the continuous case is derived. Here we will use a different technique of modulated Fourier expansions~\cite{hairer02gni}, which recently has been used  for charged-particle dynamics in a strong magnetic field~\cite{hairer20lta,hairer20afb,hairer2022large,lubich2022large,wang21eeo}, to derive the guiding centre drifts of the exact solution in such toroidal geometry. Since this technique can be extended to numerical discretization equally, the analysis of the modified Boris algorithm is also performed.

In Section~\ref{sec:setting} we formulate the equations of motion in a strongly non-uniform strong magnetic field, describe the concerned toroidal axi-symmetric  geometry and the electromagnetic field, and introduce the modified Boris scheme.
In Section~\ref{sec:main} we state the main results of this paper: Theorem~\ref{thm:exact} states the slow drift motion  over $O(\eps^{-1})$ in toroidal geometry for the continuous system and Theorem~\ref{thm:num} states the long-time accuracy of the modified Boris algorithm. Section~\ref{sec:num} presents experiments that illustrate the theoretical results. In Section~\ref{sec:proof} we give the proofs for our main results. 

\section{Setting}\label{sec:setting}
\subsection{Charged-particle dynamics in toroidal geometry}
\label{subsec:setting}
We consider the differential equation that describes  the motion of a charged particle (of unit mass and charge) in a magnetic and electric field, 
\begin{equation}\label{ode}
\ddot x =  \dot x \times B(x) + E(x),
\end{equation}
where $x(t)\in\mathbb R^3$ is the position at time $t$, $v(t)=\dot x(t)$ is the velocity, $B(x)$ is the magnetic field and $E(x)$ is the electric field.  $B$ and $E$ can be expressed via vector potential $A(x)\in\real^3$ and scalar potential $\phi(x)\in\real$ as $B(x) = \nabla \times A(x)$  and $E(x) = - \nabla \phi(x)$.  Here we are interested in the situation of a strong magnetic field
\begin{equation}\label{B-eps}
B(x) = B_\eps(x)= \frac 1\eps \, B_1(x),  \quad\ 0<\eps\ll 1,
\end{equation}
where $B_1$ is smooth and independent of the small parameter $\eps$, with $|B_1(x)|\ge 1$ for all $x$. The initial values $(x(0),\dot x(0))$ are bounded independently of $\eps$: for some constants $M_0,M_1$,
\begin{equation} \label{xv-bounds}
|x(0)| \le M_0, \quad\ |\dot x(0)| \le M_1.
\end{equation}

In this paper, we consider the so called toroidal axi-symmetric geometry which is introduced in~\cite{filbet2020asymptotics}. 
To be specific, fix a unitary vector $\e_z$,  and for  any vector $x\in\mathbb R^3$, it can be expressed as
\[
x=r(x)\,\e_r(x)+z(x)\,\e_z
\]
with $z(x)= \e_z^\top x$, $r(x)=|\e_z\times x|$, and $\e_r(x)=(x-z(x)\,\e_z)/r(x)$. It is assumed that far from the axis $\e_z$ the magnetic field is stationary, toroidal, axi-symmetric and non vanishing, that is, for some $r_0>0$, when $r(x)\geq r_0$
\begin{equation}\label{eq:B}
\e_\parallel(x)=\frac{\e_z\times x}{r(x)} \quad \text{and} \quad |B_1(x)|=b(r(x),z(x))
\end{equation}
for some function $b$. The electric field satisfies $E_\parallel(x)=0$ and $E$ is axi-symmetric when $r(x)\geq r_0$, that is,
\begin{equation}\label{eq:E}
E(x)=E_\perp(x)=E_r(r(x),z(x))\,\e_r(x)+E_z(r(x),z(x))\,\e_z.
\end{equation}
In our proofs, we assume that the functions $b$, $E_r$, $E_z$ and all their derivatives are bounded independently of $\eps$. 

\begin{figure}[h]
\centerline{
\includegraphics[scale=0.5]{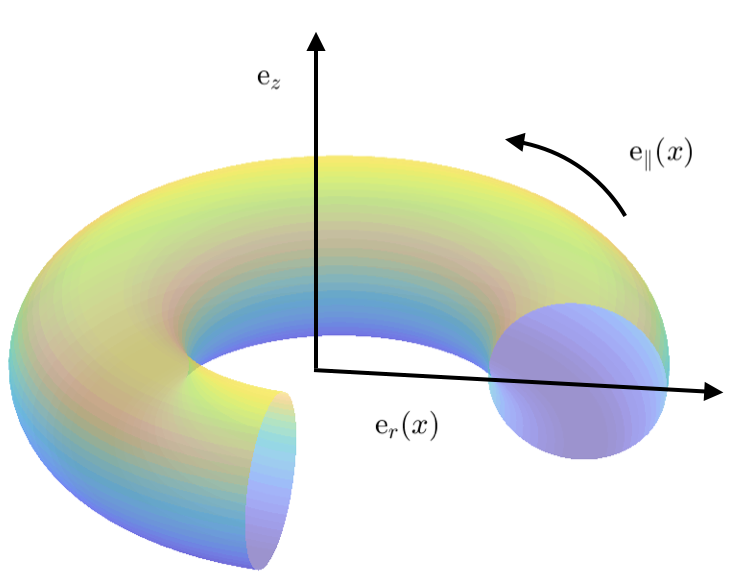}}
\caption{Toroidal geometry with $\e_r(x),\e_\parallel(x),\e_z$ the local frame  and the magnetic field along $\e_\parallel(x)$.}
\end{figure}

It is noted that $(\e_r(x), \e_\parallel(x), \e_z)$ forms the orthonormal basis and 
\[
\begin{aligned}
\e_r(x) &= \left(\frac{x_1}{r}, \frac{x_2}{r}, 0 \right)^\top\quad \e_\parallel(x)=\left(-\frac{x_2}{r}, \frac{x_1}{r}, 0\right)^\top\quad 
\e_z&=(0,0,1)^\top
\end{aligned}
\]
with $r=\sqrt{x_1^2+x_2^2}$.  The following relations  are useful in our proof
\begin{equation}\label{eq:relation}
\begin{aligned}
&\e'_r(x)=\frac{1}{r(x)}\e_\parallel(x) \e_\parallel(x)^\top, \quad \e'_\parallel(x)=-\frac{1}{r(x)} \e_r(x) \e_\parallel(x)^\top\\
&\nabla_x r(x)= \e_r(x),  \quad  B_1'(x)=\e_\parallel (\nabla_x b)^\top-\frac{b(r(x),z(x))}{r(x)}\e_r(x) \e_\parallel(x)^\top,
\end{aligned}
\end{equation}
where $'$ denotes the Jacobian of the functions considered and $\nabla_x$ is the gradient.

\subsection{Modified Boris method }
The modified Boris method proposed in \cite{xiao21smc} is used to solve the charged-particle dynamics under a strong magnetic field with large stepsizes.
The analysis of accuracy order of such method for general non-uniform strong magnetic field with stepsize $h^2\sim \eps$ until $T=O(1)$ is recently provided in~\cite{lubich2022large}. 

This algorithm has the following two-step formulation
\begin{equation}\label{mboris}
\frac{x^{n+1}-2x^n+x^{n-1}}{h^2}=v^n \times B(x^n) + E(x^n)- \mu^0\, \nabla |B|(x^n)
\end{equation}
with the initial magnetic moment 
\[\mu^0=\mu(x(0),\dot x(0))=\frac{1}{2}\frac{|\dot{x}(0)\times B(x(0))|^2}{|B(x(0))|^3}.
\] The velocity is computed as
\begin{equation} \label{vn}
v^n = \frac{x^{n+1}-x^{n-1}}{2h}.
\end{equation}
The modified Boris method starts from modified initial values 
\begin{equation}\label{mod-init}
x^0 = x(0), \quad v^0 = P_\parallel(x^0) \,\dot x(0),
\end{equation}
where $P_\parallel(x^0)$ is the orthogonal projection onto the span of $B(x^0)$. This means the component of the initial velocity orthogonal to the magnetic field is filtered out, i.e., $v^0_\perp=P_\perp(x^0)v^0=0$ with $P_\perp(x^0)=I-P_\parallel(x^0)$.

We note that the modified Boris method is identical to the standard Boris integrator for the modified electric field $E_\mathrm{mod}(x) = E(x)- \mu^0 \,\nabla |B|(x) =-\nabla (\phi + \mu^0 |B|)(x)$ and can be implemented as the common one-step formulation of the Boris algorithm~\cite{boris70rps}.

\section{Main results}\label{sec:main}
Introducing  $\tilde{r}(t), \tilde{z}(t), \tilde{v}(t)$ such that they are the solutions of the following initial-value problem for the slow differential equations
\begin{equation}\label{eq:limit}
\begin{aligned}
\frac{\d \tilde r}{\d t}&=-\eps\frac{E_z(\tilde{r},\tilde{z})}{b(\tilde{r},\tilde{z})}+\eps\frac{\mu^0}{b(\tilde{r},\tilde{z})}\partial_z b(\tilde{r},\tilde{z}), \quad \tilde{r}(0)=r(x(0))\\
\frac{\d \tilde z}{\d t}&=\eps\frac{\tilde{v}^2}{\tilde{r} b(\tilde{r},\tilde{z})}+\eps\frac{E_r(\tilde{r},\tilde{z})}{b(\tilde{r},\tilde{z})}-\eps\frac{\mu^0}{b(\tilde{r},\tilde{z})}\partial_r b(\tilde{r}, \tilde{z}), \quad \tilde{z}(0)=z(x(0))\\
\frac{\d \tilde v}{\d t}&=\eps\frac{\tilde{v}}{\tilde{r}}\left(\frac{E_z(\tilde{r}, \tilde{z})}{b(\tilde{r}, \tilde{z})}-\frac{\mu^0}{b(\tilde{r}, \tilde{z})}\partial_z b(\tilde{r}, \tilde{z})\right), \quad \tilde{v}(0)=\e_\parallel(x(0))^\top \dot{x}(0),
\end{aligned}
\end{equation}
we then have the following results.

\begin{theorem}[Drift motion of the exact solution]\label{thm:exact}
Let $x(t)=r(x(t))\,\e_r(x(t))+z(x(t))\,\e_z$ be a solution of \eqref{ode}--\eqref{xv-bounds} with \eqref{eq:B} and \eqref{eq:E}, which stays in a compact set $K$ for $0\leq t\leq c\eps^{-1}$ (with $K$ and $c$ independent of $\eps$) and $v_\parallel(t)= \e_\parallel(x(t))^\top \dot{x}(t)$ be the parallel velocity.  Denote $r(t)=r(x(t))$ and $z(t)=z(x(t))$, then we have
\[
|r(t)-\tilde{r}(t)|\leq C\eps, \  |z(t)-\tilde{z}(t)|\leq C\eps, \  |v_\parallel(t)-\tilde{v}(t)|\leq C\eps, \quad 0\leq t \leq c/\eps.
\]
The constant $C$ is independent of $\eps$ and $t$ with $0\leq t\leq c/\eps$, but depends on $c$ and on bounds of derivatives of $B_1$ and $E$ on the compact set $K$.
\end{theorem}
\begin{remark} A similar result is given by Proposition 5.2 in~\cite{filbet2020asymptotics}. Here we provide a different proof of the modulated Fourier expansions, which can be extended to the analysis of numerical methods and enables us to obtain the following result.
\end{remark}
For the numerical approximation, the nondegeneracy condition is needed as in~\cite{lubich2022large}:
\begin{align}\label{ass-nondeg}
&\text{For $(x,v)$ along the numerical trajectory, the linear maps} \nonumber
\\[1mm]
& \text{$L_{x,v}:P_\perp(x)\mathbb{R}^3 \to P_\perp(x)\mathbb{R}^3, \quad z \mapsto z + \tfrac14 h^2\, P_\perp(x)\bigl(v \times B'(x)z\bigr)$}
\\[1mm]
&\text{have an inverse that is bounded independently of $(x,v)$ and of 
} \nonumber
\\[-1mm]
&\text{$h$ and $\eps$ with $h^2/\eps\le C_*$.} \nonumber
\end{align}
This determines an upper bound $C_*$ on the ratio $h^2/\eps$.

\begin{theorem}[Drift approximation by the numerical solution]\label{thm:num}
Consider applying the modified Boris method to \eqref{ode}--\eqref{xv-bounds} with \eqref{eq:B} and \eqref{eq:E} and with modified initial values \eqref{mod-init} using a step size $h$ with $h^2\sim \eps$, i.e.,
\[
c_* \eps \le h^2 \le C_* \eps
\]
for some positive constants $c_*$ and $C_*$. Under the nondegeneracy condition \eqref{ass-nondeg} and provided that the numerical solution $x^n=r(x^n)\,\e_r(x^n)+z(x^n)\,\e_z$ stays in a compact set $K$ for $0\leq nh\leq c\eps^{-1}$ (with $K$ and $c$ independent of $\eps$ and $h$). $v_\parallel^n=\e_\parallel(x^n)^\top v^n$ denotes the parallel component  of the numerical velocity $v^n$. Then
\[
|r(x^n)-\tilde{r}(t_n)|\leq Ch^2, \ |z(x^n)-\tilde{z}(t_n)|\leq Ch^2, \  |v_\parallel^n-\tilde{v}(t_n)|\leq Ch^2,  \quad 0\leq t_n=nh\leq c/\eps.
\]
The constant $C$ is independent of $\eps$ and $h$ and $n$ with $0\leq nh\leq c/\eps$, but depends on on $c$, on bounds of derivatives of $B_1$ and $E$ on the compact set $K$ and on $c_*$ and $C_*$. 

\end{theorem}

\section{Numerical experiments}\label{sec:num}
\begin{figure}[htbp]
\centerline{
\includegraphics[scale=0.5]{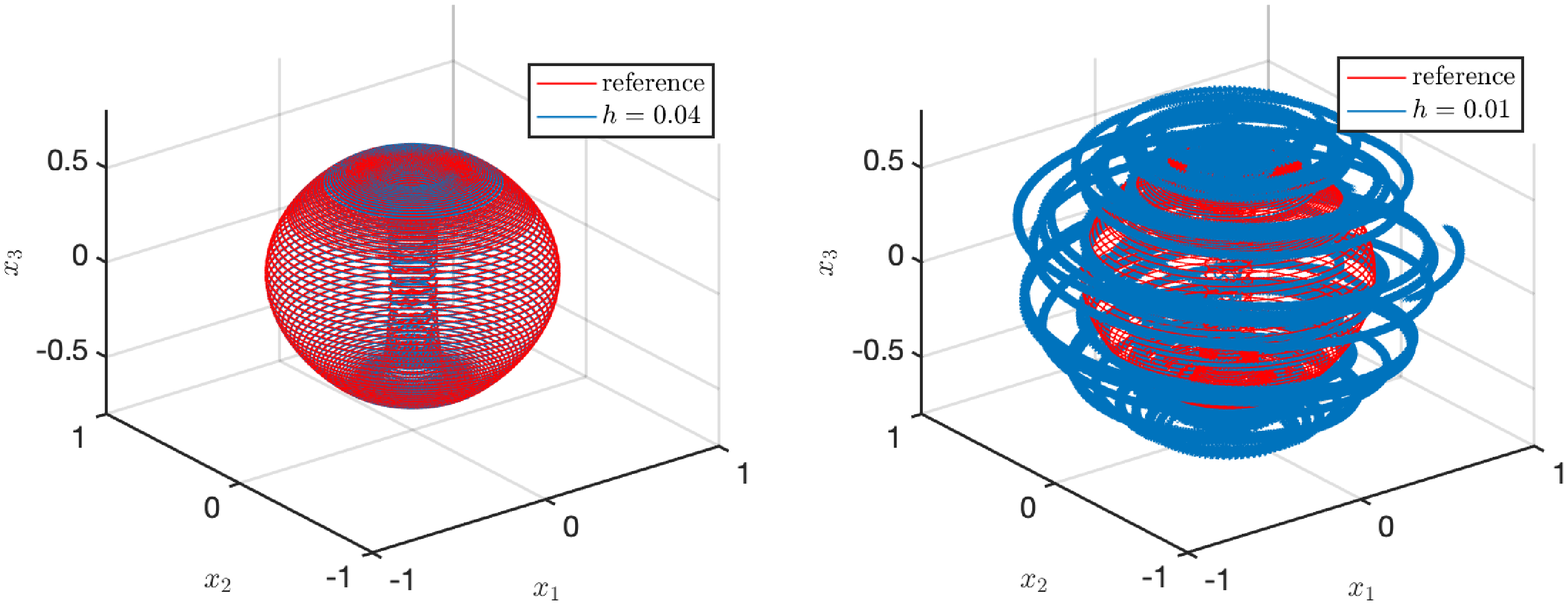}}
\caption{Particle trajectories for $t\leq 1/\eps$ with $\eps=10^{-3}$ as computed by the modified Boris with $h=0.04$ (left) and by the Boris method with $h=0.01$ (right). }\label{fig:orbit1}
\end{figure}

\begin{figure}[htbp]
\centerline{
\includegraphics[scale=0.5]{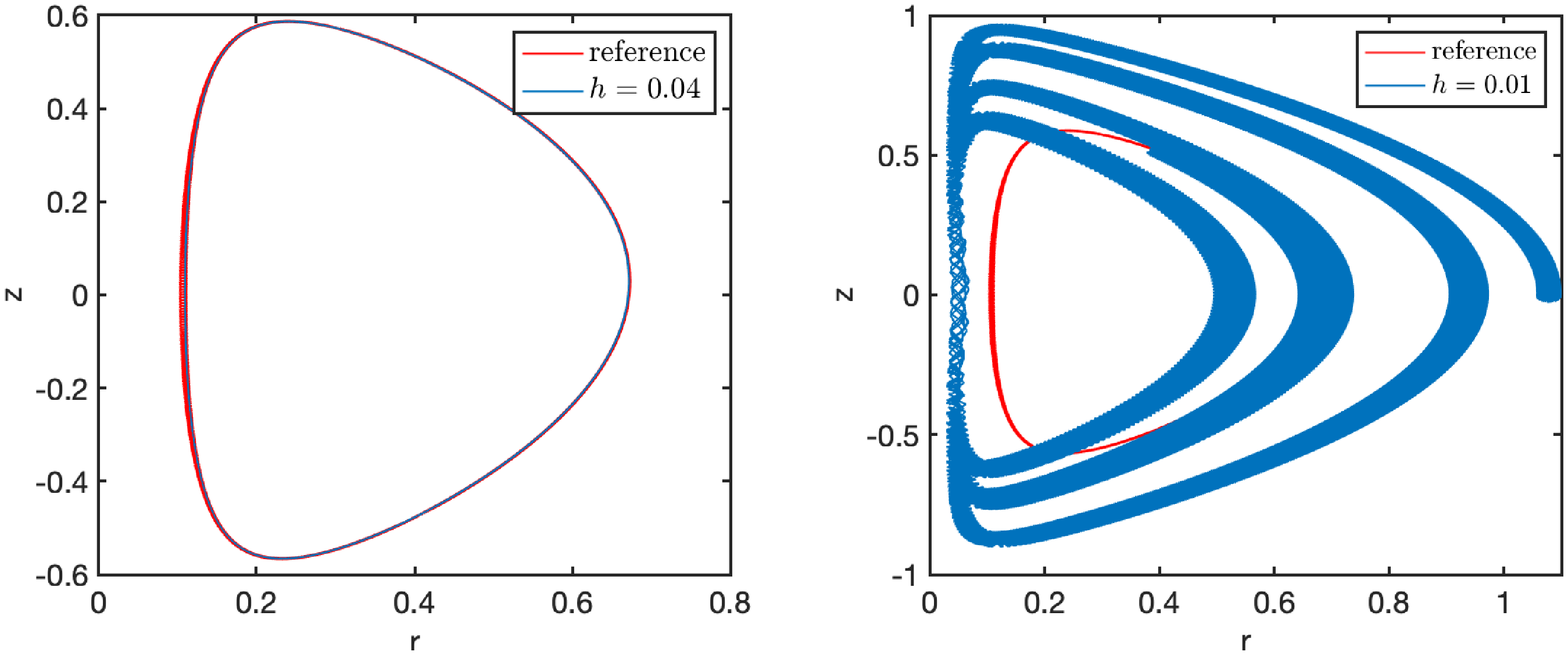}}
\caption{Particle trajectories for $t\leq 1/\eps$ with $\eps=10^{-3}$ projected onto the $r-z$ plane as computed by the modified Boris with $h=0.04$ (left) and by the Boris method with $h=0.01$ (right). }\label{fig:orbit2}
\end{figure}

\begin{figure}[htbp]
\centerline{
\includegraphics[scale=0.5]{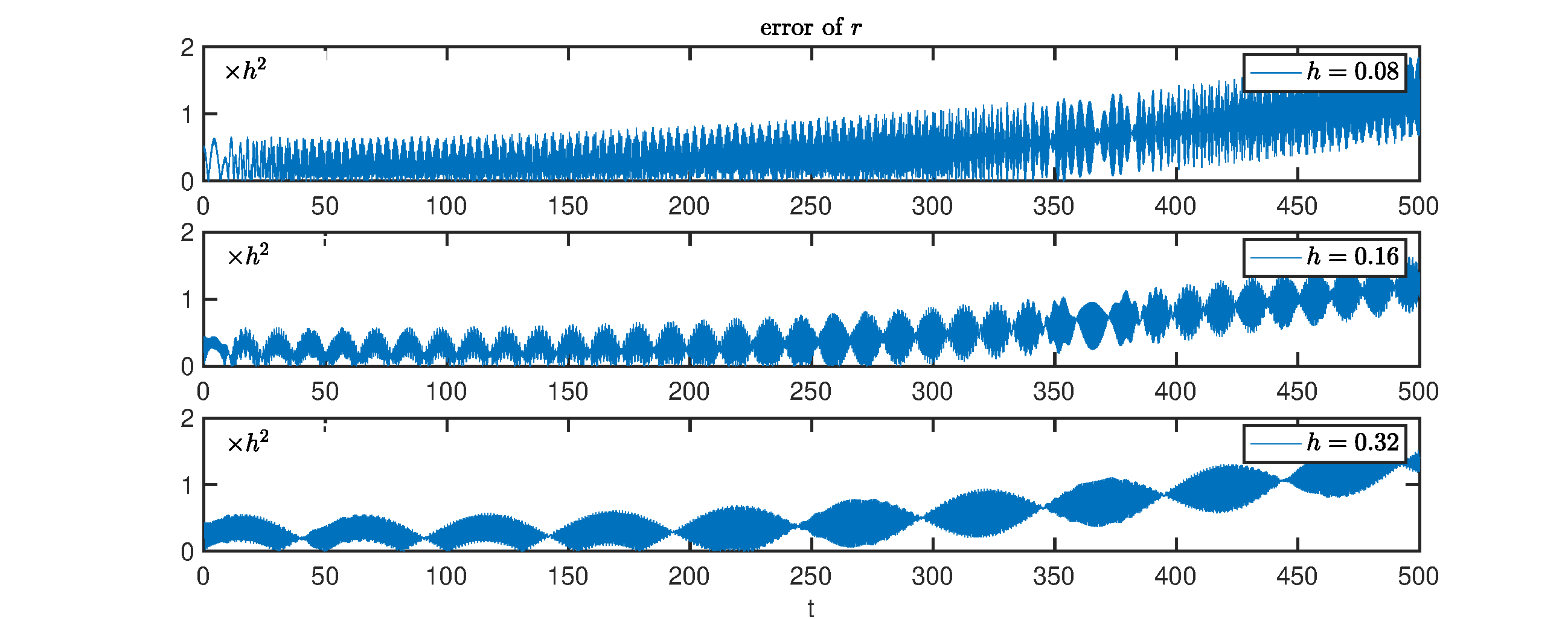}}
\centerline{
\includegraphics[scale=0.5]{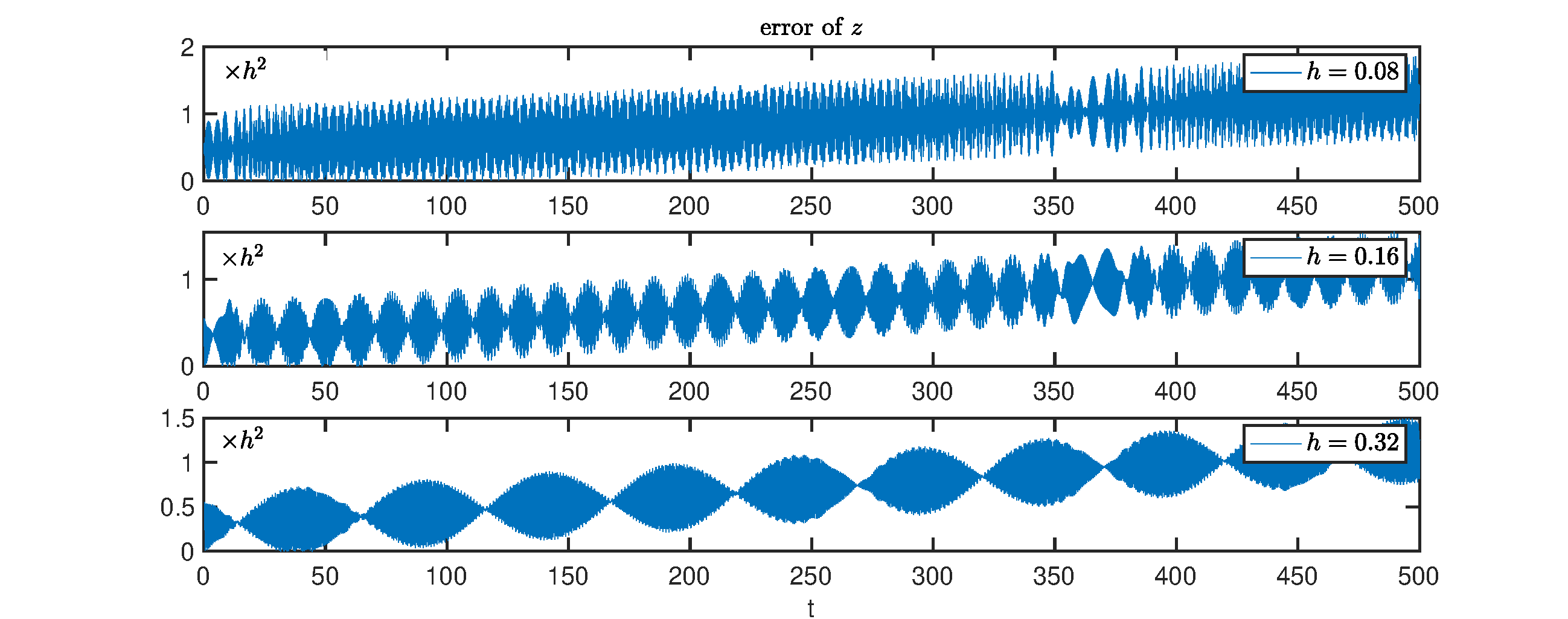}}
\centerline{
\includegraphics[scale=0.5]{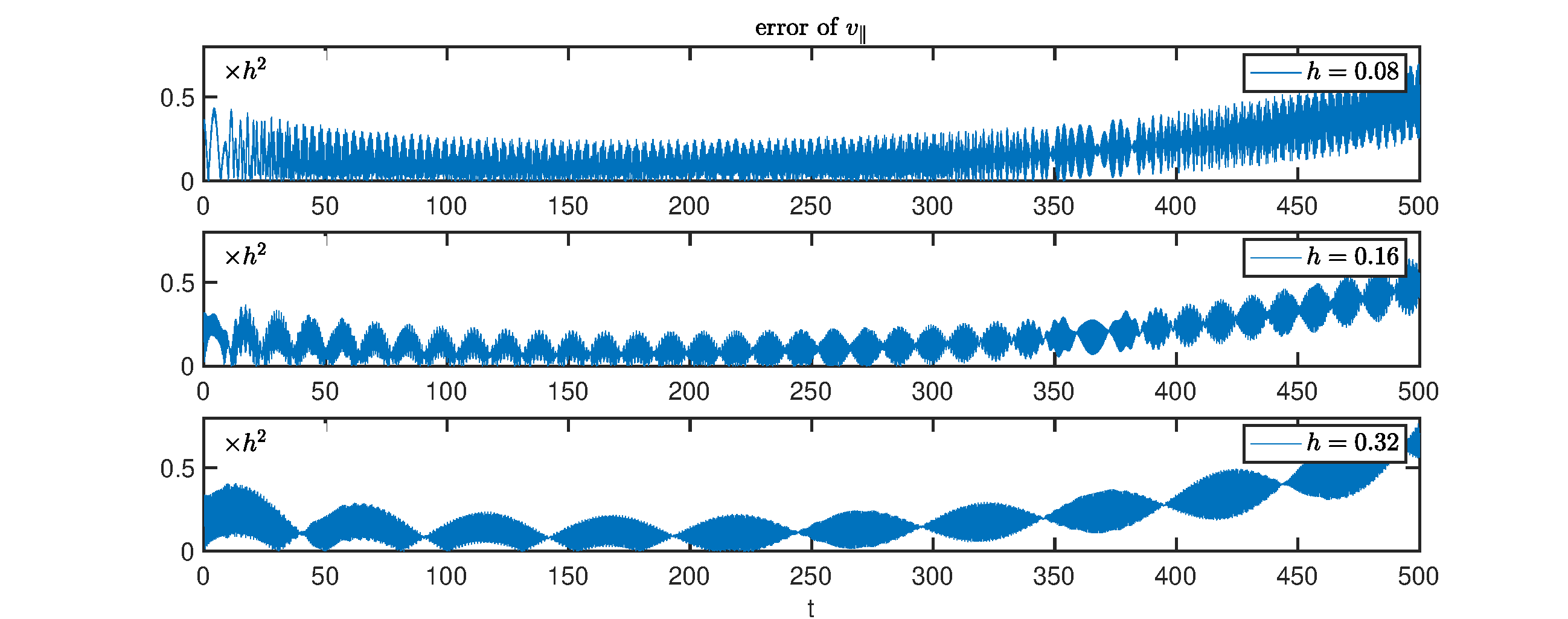}}
\caption{Absolute errors $|r(x^n)-r^{\text{ref}}|$, $|z(x^n)-z^{\text{ref}}|$ and $|v_\parallel^n-v_\parallel^{\text{ref}}|$ as functions of time, along the numerical solution of the modified Boris algorithm with $\eps=10^{-3}$ and three different $h$.}\label{fig:order1}
\end{figure}

\begin{figure}[htbp]
\centerline{
\includegraphics[scale=0.5]{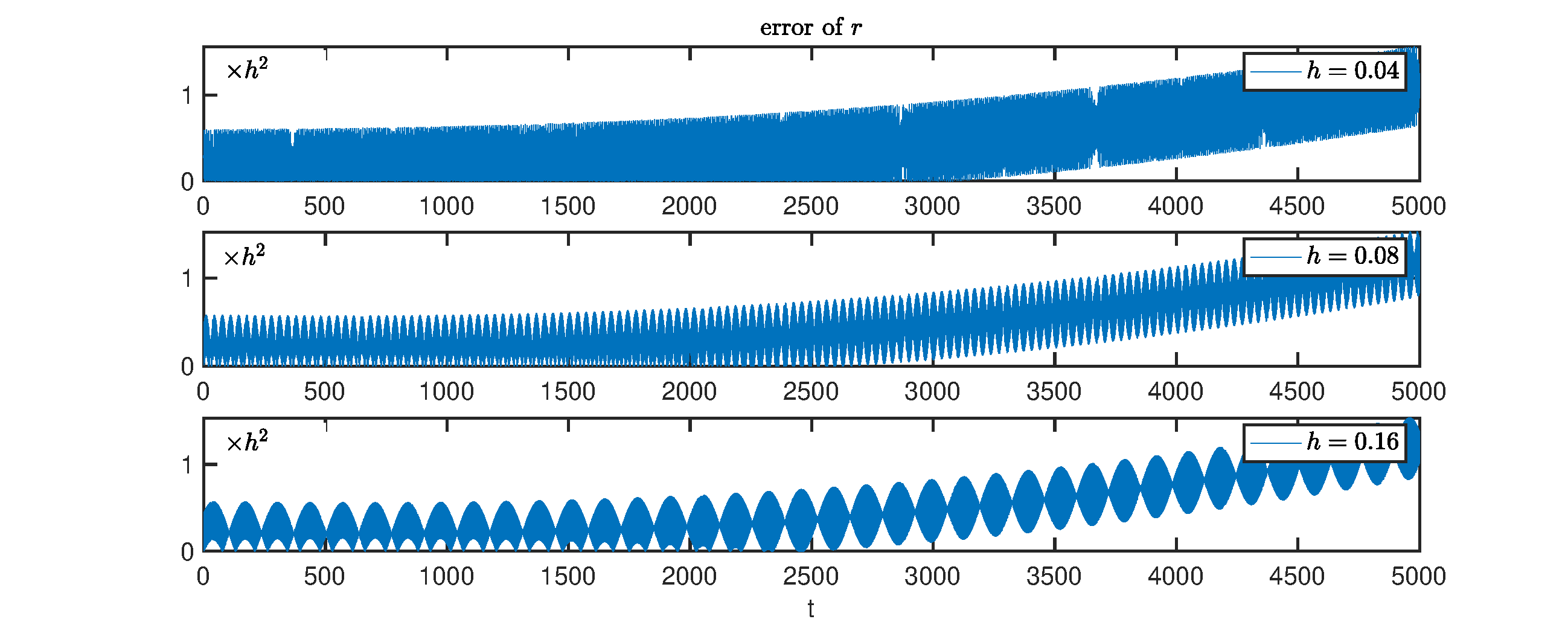}}
\centerline{
\includegraphics[scale=0.5]{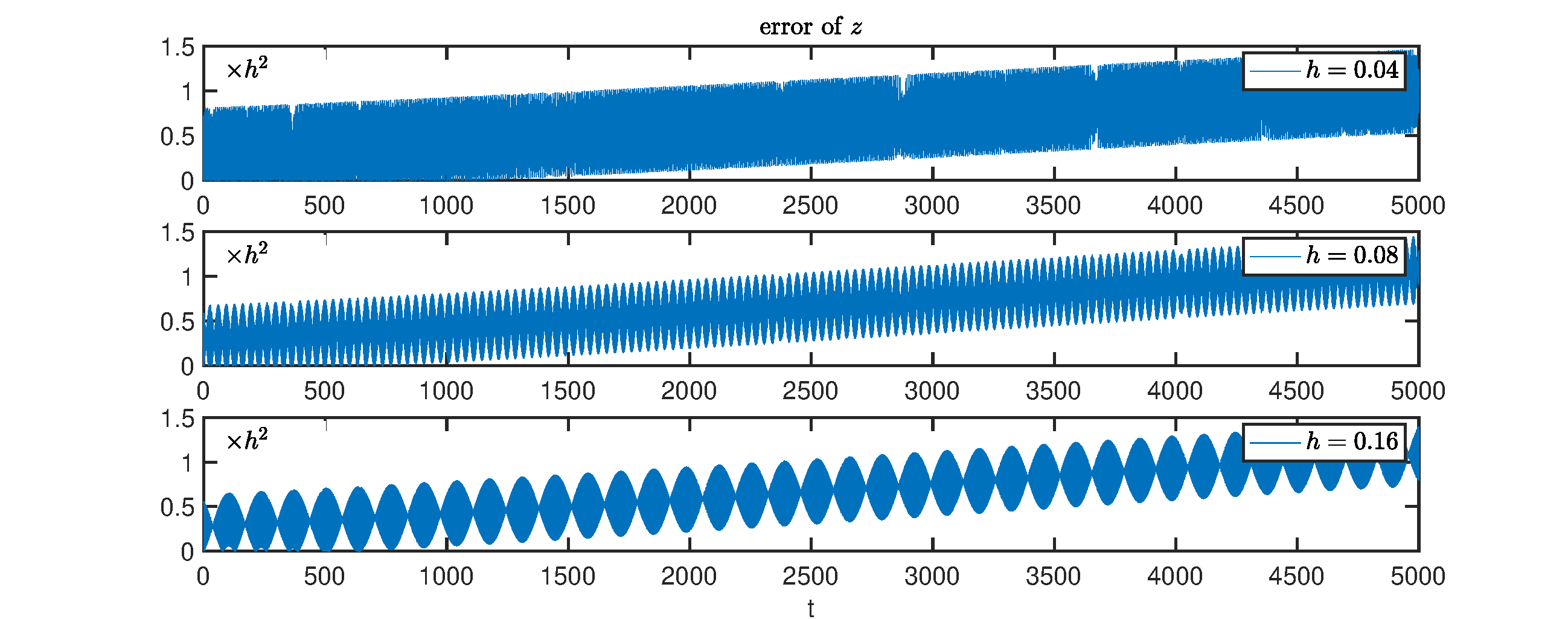}}
\centerline{
\includegraphics[scale=0.5]{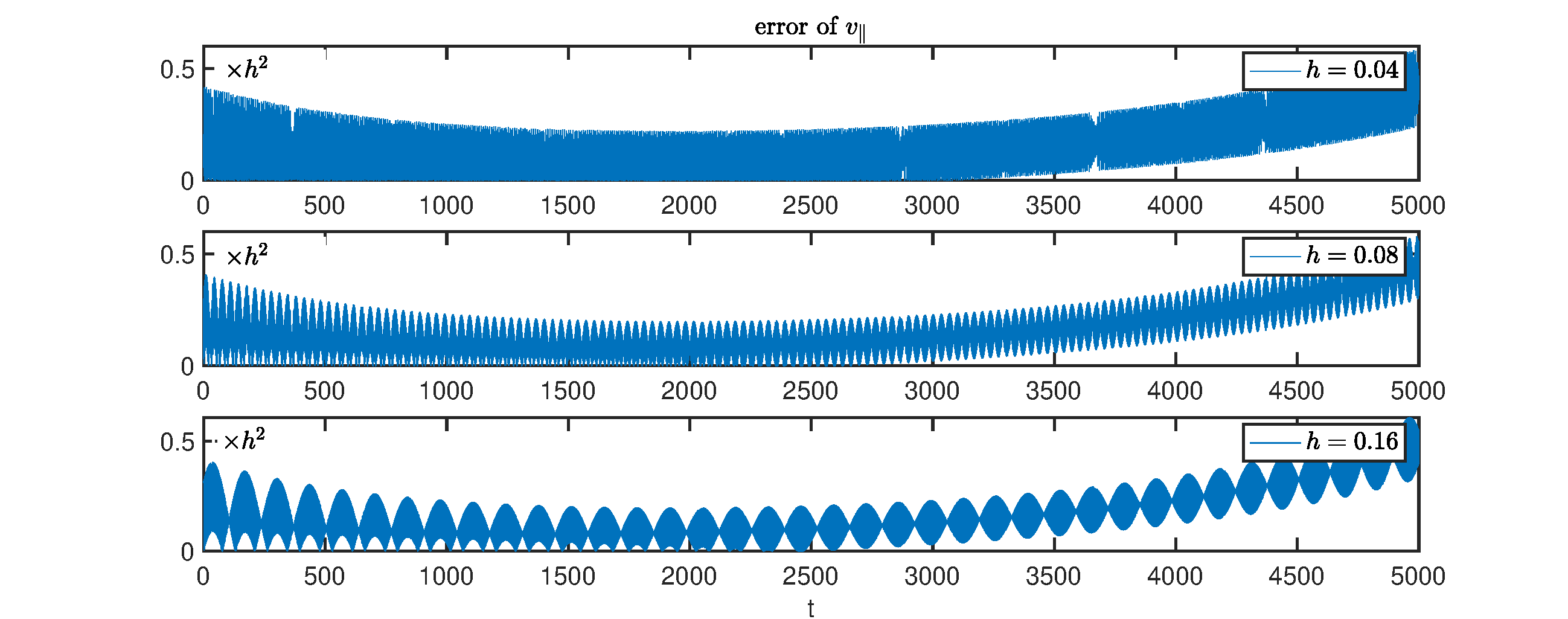}}
\caption{Absolute errors $|r(x^n)-r^{\text{ref}}|$, $|z(x^n)-z^{\text{ref}}|$ and $|v_\parallel^n-v_\parallel^{\text{ref}}|$ as functions of time, along the numerical solution of the modified Boris algorithm with $\eps=10^{-4}$ and three different $h$.}\label{fig:order2}
\end{figure}

To illustrate the statement of the preceding section we consider the following electromagnetic fields
\[
\begin{aligned}
E(x)&=0.1 z(x) \, \e_r(x)+0.1 r(x)\, \e_z=0.1\left (\frac{x_1 x_3}{r}, \,  \frac{x_2 x_3}{r},\, r\right)^\top, \\
B(x)&=\frac{r(x)+z^2(x)}{\eps}\e_\parallel(x)=\frac{r+x^2_3}{\eps}\left(-\frac{x_2}{r},\, \frac{x_1}{r}, \, 0\right)^\top
 \end{aligned}
\]
with $r=\sqrt{x_1^2+x_2^2}$. The initial values are chosen as
\[
x(0)=(1/3,1/4,1/2)^\top, \quad  \dot{x}(0)=(2/5, 2/3,1)^\top.
\]
Figure~\ref{fig:orbit1} shows the trajectories computed by the standard Boris and modified Boris with final time $T=1/\eps$. It is observed that the modified Boris method can give the correct trajectories even with very large time step size $h=40\eps$. For the standard Boris method, the drift motions are totally wrong with large time step size. The projection of the computed particle trajectory onto the $(r,z)$ plane is given in Figure \ref{fig:orbit2}. Figure \ref{fig:order1} shows the absolute errors of $r$, $z$, and $v_\parallel$ along the numerical solution of the modified Boris algorithm with $\eps=10^{-3}$ and $T=0.5/\eps$, which is observed to be size of $O(h^2)$ in agreement  with our theoretical results. Figure \ref{fig:order2} shows the similar results for $\eps=10^{-4}$. All the reference solutions are obtained by using standard Boris with small time step size $h=0.05\eps$.

\section{Proof of main results}~\label{sec:proof}
The theorems will be proved mainly based on the modulated Fourier expansions for the exact and numerical solutions given in~\cite{lubich2022large}. In this section, we will write the guiding centre equations in the toroidal geometry and express all the $O(\eps)$ terms explicitly.

Following~\cite{hairer20lta}, we diagonalize the linear map $v\mapsto v\times B(x)$ and denote the eigenvalues as $\lambda_1=\iu|B(x)|$, $\lambda_0=0$, and $\lambda_{-1}=-\iu|B(x)|$. The corresponding normalised eigenvectors are denoted by $\nu_1(x), \, \nu_0(x), \, \nu_{-1}(x)$ and the  orthogonal projections onto the eigenspaces are denoted by $P_j(x)=\nu_j(x)\nu_j(x)^*$. It is noted that $P_\parallel(x)=P_0(x)$ and $P_\perp(x)=I-P_\parallel(x)=P_1(x)+P_{-1}(x)$.
\subsection{Proof of Theorem~\ref{thm:exact}}
The proof is structured into three parts (a)-(c).

(a) The equation of guiding centre motion in cartesian coordinate.

According to Theorem 4.1 of \cite{hairer20lta}, it is known that the solution of \eqref{ode}--\eqref{B-eps} can be written as
\[
x(t)=\sum_{|k|\leq N-1}y^k(t)\,\e^{{\iu}k\varphi(t)/\eps} + R_N(t), \qquad 0\le t \le c\eps,
\]
where the phase function satisfies $\dot{\varphi}(t)=|B_1(y^0(t))|$. The coefficient functions $y^k(t)$ together with their derivatives (up to order $N$) are bounded as
\[
y^k=O(\eps^{|k|}) \quad \forall |k|\leq N-1
\]
and further satisfy 
\[
\dot y^0 \times B_1(z^0)=O(\eps), \quad y^k_j  =O(\eps^2) \quad \text{for} \quad |k|=1, j\neq k,
\]
where $y^k(t)=y^k_1(t)+y^k_0(t)+y^k_{-1}(t)$ with $y^k_j=P_j(y^0)y^k$.
The remainder term and its derivative are bounded by 
\[
R_N(t)=O(\eps^N),\quad \dot{R}_N(t)=O(\eps^{N-1}).
\] 

Similar to Theorem 4.1 of \cite{lubich2022large}, we can divide the interval $[0, c/\eps]$ into small intervals of length $O(\eps)$ and on each subinterval we consider the above modulated Fourier expansion, which means $x(t)$ can be written as the modulated Fourier expansion for longer time intervals 
\begin{equation}\label{eq:mfe}
x(t)=\sum_{|k|\leq N-1}y^k(t)\,\e^{{\iu}k\varphi(t)/\eps} + R_N(t), \qquad 0\le t \le \frac{c}{\eps},
\end{equation}
where $y^k(t)$ are piecewise continuous with jumps of size $O(\eps^{N})$ at integral multiples of $\eps$ and are smooth elsewhere. The sizes of the coefficients and remainder term are the same as above.

Inserting \eqref{eq:mfe} into the continuous system and comparing the coefficients of $\e^{\iu k\varphi(t)/\eps}$ yield the differential equations for $y^k(t)$. For $k=0$ and  $k=\pm 1$, we have
\begin{equation}\label{eq:eqs}
\begin{aligned}
\ddot{y}^0=\dot{y}^0\times B(y^0)+E(y^0)+\underbrace{2\mathrm{Re}\left(\iu|B|y^1\times B'(y^0)y^{-1}\right)}_{=:I}+\underbrace{2\mathrm{Re}\left(\dot{y}_1^1\times B'(y^0)y^{-1}_{-1}\right)}_{=:II}+O(\eps^2),\\
\pm 2\iu\frac{\dot{\varphi}}{\varepsilon}\dot{y}^{\pm1}+\left(\pm\iu\frac{\ddot{\varphi}}{\eps}-\frac{\dot{\varphi}^2}{\eps^2}\right)y^{\pm1}=\left(\dot{y}^{\pm1}\pm \iu\frac{\dot{\varphi}}{\eps}y^{\pm1}\right)\times B(y^0)+\dot{y}^0\times B'(y^0)y^{\pm1}_{\pm1}+O(\eps).
\end{aligned}
\end{equation}
From the first equation of \eqref{eq:eqs}, it is straightforward to get  several slow drifts for $P_\perp\dot{y}^0$ (see Remark 4.3 of \cite{lubich2022large}) and the guiding centre motion of $y^0(t)$ satisfies
\begin{equation}\label{eq:y}
\dot{y}^0=P_{\parallel}\dot{y}^0+\frac{1}{|B|}P_{\parallel}\dot{y}^0\times \frac{\d \e_\parallel}{\d t}+\frac{1}{|B|^2}\left(E-\mu^0\,\nabla |B|\right)\times B+O(\eps^2),
\end{equation}
with $B, \nabla |B|, P_\parallel=\e_\parallel \e_\parallel^\top$ and $E$ evaluated at the guiding centre $y^0$. The initial value of $y^0$ is 
\[
y^0(0)= x(0)+ \frac{\dot{x}(0)\times B(x(0))}{|B(x(0))|^2}+O(\varepsilon^2).\\
\]

(b) The equations in toroidal geometry.

In the toroidal  geometry, $y^0(t)$ can be written as 
\[
y^0=r(y^0)\e_r(y^0)+z(y^0)\e_z=:r^0\e_r(y^0)+z^0 \e_z.
\]

\noindent
--- Multiplying  \eqref{eq:y} with  $\e_r^\top=\e_r(y^0)^\top$ gives
\begin{equation}\label{eq:r}
\e_r^\top\dot{y}^0 =\frac{\eps v^0_\parallel}{b} \e_r^\top\left(\e_\parallel\times\frac{\d \e_\parallel}{\d t}\right)+\frac{\eps}{b} \e_r^\top \left((E-\mu^0\,\nabla b)\times \e_\parallel\right) +O(\eps^2),
\end{equation}
where $\e_r, \e_\parallel$ are evaluated at $y^0$  and
\begin{equation}\label{eq:vv}
v^0_\parallel:= \e_\parallel^\top\dot{y}^0.
\end{equation}
From \eqref{eq:relation}, it is known that
\begin{equation}\label{eq:deriv}
\begin{aligned}
\dot{\e}_r(y^0)&=\frac{v^0_\parallel}{r^0}\e_\parallel(y^0), \quad \dot{\e}_\parallel(y^0)=-\frac{v^0_\parallel}{r(y^0)} \e_r(y^0),\\ \ddot{\e}_\parallel(y^0)&=-\frac{\d}{\d t}\left(\frac{v^0_\parallel}{r^0}\right)\e_r(y^0)-\left(\frac{v^0_\parallel}{r^0}\right)^2\e_\parallel(y^0).
\end{aligned}
\end{equation}
Then the left hand side of \eqref{eq:r} can be  expressed as
\[
\e_r^\top \dot{y}^0 = \frac{\d}{\d t}( \e_r^\top y^0 ) -  \dot{\e}_r^\top y^0=\frac{\d r^0}{\d t},
\]
and the first term on the right hand side of \eqref{eq:r} vanishes since
\[
\e_r^\top \left( \e_\parallel\times \frac{\d \e_\parallel}{\d t} \right)= -\frac{v^0_\parallel}{r} \e_r^\top( \e_\parallel\times \e_r)=0.
\]
Using the fact that $\e_r\times\e_\parallel=\e_z, \, \e_z\times \e_\parallel=-\e_r$, $E=E_r\e_r+E_z\e_z$ and $\nabla b=\partial_r b \, \e_r+\partial_z b \, \e_z$, we obtain
\[
\e_r^\top \left((E-\mu^0\,\nabla b(r^0, z^0))\times \e_\parallel\right) = -E_z+\mu^0\partial_z b.
\]
Thus \eqref{eq:r} is equivalent to
\[
\frac{\d r^0}{\d t}=\frac{\eps}{b}( -E_z+\mu^0\partial_z b)+O(\eps^2),
\]
where the functions $E_z, b, \partial_z b$ are evaluated at $(r^0,z^0)$.
The initial value of $r^0$ can be expressed as
\[
r^0(0)=\e_r(y(0))^\top y^0(0)=\e_r(x(0))^\top x(0) + O(\eps)=r(x(0))+O(\eps).
\]

\noindent
--- Multiplying \eqref{eq:y} with $\e_z^\top$ gives
\begin{equation}\label{eq:z}
\e_z^\top \dot{y}^0 =\frac{\eps v^0_\parallel}{b(r^0, z^0)} \e_z^\top \left(\e_\parallel\times \frac{\d \e_\parallel}{\d t}\right)+\frac{\eps}{b(r^0, z^0)} \e_z^\top \left((E-\mu^0\,\nabla b(r^0, z^0))\times \e_\parallel\right) +O(\eps^2)
\end{equation}
Similarly, we have
\[
\e_z^\top \dot{y}^0 = \frac{\d}{\d t}(\e_z^\top y^0)=\frac{\d z^0}{\d t},
\]
\[
\e_z^\top\left( \e_\parallel\times \frac{\d \e_\parallel}{\d t}\right) = -\frac{v^0_\parallel}{r^0}\e_z^\top( \e_\parallel\times \e_r)=\frac{v^0_\parallel}{r^0},
\]
and
\[
\e_z^\top \left((E-\mu^0\,\nabla b)\times \e_\parallel\right) = E_r-\mu^0\partial_r b,
\]
then \eqref{eq:z} can be expressed as
\[
\frac{\d z^0}{\d t}=\eps\frac{(v^0_\parallel)^2}{r^0 b}+\eps\frac{E_r}{b}-\eps\frac{\mu^0}{b}\partial_r b+O(\eps^2),
\]
where the functions $E_r, b, \partial_r b$ are evaluated at $(r^0,z^0)$.
The initial value of $z^0$ is
\[
z^0(0)=\e_z^\top y^0(0)=\e_z^\top x(0) + O(\eps)=z(x(0))+O(\eps).
\]

\noindent
--- By the definition of \eqref{eq:vv} we can derive the equation for $v^0_\parallel$ 
\begin{equation}\label{eq:dvv}
\frac{\d}{\d t}v^0_\parallel=\frac{\d}{\d t} (\e_\parallel^\top \dot{y}^0)= \dot{\e}_\parallel^\top \dot{y}^0+ \e_\parallel^\top \ddot{y}^0.
\end{equation}
The first term on the right hand side is 
\begin{equation}\label{eq:form}
\dot{\e}_\parallel^\top \dot{y}^0=-\frac{v^0_\parallel}{r^0}\frac{\d r^0}{\d t}
\end{equation}
using \eqref{eq:deriv}.
In the following we will show that the second term $\e_\parallel^\top \ddot{y}^0$ is of size $O(\eps^2)$.

From the expression of $B'$ given in \eqref{eq:relation}, it is known that $ B'(y^0)y^{k}_{\pm1}$ is parallel to $\e_\parallel$, thus $ \e_\parallel^\top II=0$ and
\[
 \e_\parallel^\top I=  \e_\parallel^\top \left(2 \, \mathrm{Re}(\iu|B|y_1^1\times B'(y^0)y_0^{-1}\right)+O(\eps^2).
\]
The algebraic equation of $y_0^{\pm1}$ can be derived by applying $P_\parallel(y^0)$ to the second equation of \eqref{eq:eqs} 
\[
\pm 2\iu\frac{\dot{\varphi}}{\varepsilon}P_\parallel\dot{y}^{\pm1}+\left(\pm\iu\frac{\ddot{\varphi}}{\eps}-\frac{\dot{\varphi}^2}{\eps^2}\right)y_0^{\pm1}=P_\parallel\left(\dot{y}^0\times B'(y^0)y^{\pm1}_{\pm1}\right)+O(\eps).
\]
The dominant term is $-\dot{\varphi}^2/\eps^2 y_0^{\pm1}$ and the right hand side $P_\parallel\left(\dot{y}^0\times B'(y^0)y^{\pm1}_{\pm1}\right)=0$ using that $B'(y^0)y^{k}_{\pm1}$ is parallel to $\e_\parallel$. 
Hence we obtain the following relation for $y_0^{\pm1}$
\[
\begin{aligned}
y_0^{\pm1}&=\pm 2\iu \frac{\eps}{\dot{\varphi}}P_\parallel\dot{y}^{\pm1}+O(\eps^3)\\
&=\pm 2\iu \frac{\eps}{\dot{\varphi}}\dot{y}^{\pm1}_0 \mp 2\iu \frac{\eps}{\dot{\varphi}}\dot{P}_\parallel y^{\pm1}_{\pm1}+O(\eps^3).
\end{aligned}
\]
By differential and substitution the first term on the right hand side of above equation can be removed. Using \eqref{eq:deriv}, we have
\[
\begin{aligned}
y_0^{\pm1}&=\mp 2\iu \frac{\eps}{\dot{\varphi}}(\dot{\e}_\parallel \e^\top_\parallel+\e_\parallel \dot{\e}^\top_\parallel) y^{\pm1}_{\pm1}+O(\eps^3)\\
&=\pm 2\iu \frac{\eps}{\dot{\varphi}}\frac{v^0_\parallel}{r^0}({\e}^\top_r y^{\pm1}_{\pm1})\e_\parallel +O(\eps^3).\\
\end{aligned}
\]
Denoting $y^{1}_{1}=\zeta\,\nu_1$, $y^{-1}_{-1}=\bar{\zeta}\,\nu_{-1}$ with $\nu_{\pm1}=(\e_z \pm\iu \,\e_r(y^0))/{\sqrt 2}$ and substituting it to the above equation give  $y^{1}_0=\eta \, \e_\parallel+O(\eps^3)$ and $y^{-1}_0=\bar{\eta}\,\e_\parallel+O(\eps^3)$  with $\eta=-\eps(\sqrt{2}{ v^0_\parallel}/{\dot{\varphi}r^0}) \zeta$. Then we have
\[
\begin{aligned}
2\,\mathrm{Re}\left(\iu|B|y_1^1\times B'(y^0)y_0^{-1}\right)&=\sqrt{2}\,\mathrm{Re}\left(\iu|B| \zeta\bar{\eta}(\e_z\times B'(y^0)\e_\parallel+\iu \,\e_r\times B'(y^0)\e_\parallel)\right)+O(\eps^2)\\
&=-\sqrt{2}\,|B| \zeta\bar{\eta} \ \e_r\times B'(y^0)\e_\parallel +O(\eps^2).
\end{aligned}
\]
From the expression of $B'$ in \eqref{eq:relation}, we konw that $B'(y^0)\e_\parallel$ is the combination of $\e_\parallel$ and $\e_r$ and thus $\e_\parallel^\top( \e_r\times B'(y^0)\e_\parallel)=0$. This means
\[
 \e_\parallel^\top \ddot{y}^0=  \e_\parallel^\top I +O(\eps^2)=O(\eps^2),
\]
and \eqref{eq:dvv} is equivalent to  
\[
\frac{\d}{\d t}v^0_\parallel=-\frac{v^0_\parallel}{r^0}\frac{\d r^0}{\d t}+ O(\eps^2).
\]
 The initial value of $v^0_\parallel$ is
\[
v^0_\parallel(0)=\e_\parallel(y^0(0))^\top \dot{y}^0(0)= \e_\parallel(x(0))^\top \dot{x}(0) + O(\eps).
\]

(c) From short to long time intervals

Denoting by $y^{0,[n]}, r^{0,[n]}, z^{0,[n]}, v^{0,[n]}_\parallel$ the functions  $y^{0}, r^0, z^0, v^0_\parallel$  on time interval $n\eps\leq t\leq (n+1)\eps$, from (b), it is known that  these coefficients satisfy the following equations
\begin{equation}\label{eq:gc}
\begin{aligned}
\frac{\d r^{0,[n]}}{\d t}&=-\eps\frac{E_z}{b}+\eps\frac{\mu^0}{b}\partial_z b+O(\eps^2), \quad r^{0,[n]}(n\eps)=r(x(n\eps))+O(\eps)\\
\frac{\d z^{0,[n]}}{\d t}&=\eps\frac{(v^{0,[n]}_\parallel)^2}{b r^{0,[n]}}+\eps\frac{E_r}{b}-\eps\frac{\mu^0}{b}\partial_r b+O(\eps^2), \quad 
z^{0,[n]}(n\eps)=z(x(n\eps))+O(\eps)\\
\frac{\d v^{0,[n]}_\parallel}{\d t}&=\eps\frac{v^{0,[n]}_\parallel}{r^{0,[n]}}\left(\frac{E_z}{b}-\frac{\mu^0}{b}\partial_z b\right)+O(\eps^2), \quad v^{0,[n]}_\parallel(n\eps)=\langle \dot{x}(n\eps), \e_\parallel(x(n\eps)) \rangle + O(\eps),
\end{aligned}
\end{equation}
with $E_r, E_z, b, \partial_r b, \partial_z b$ evaluated at $(r^{0,[n]},z^{0,[n]})$.
From equation \eqref{eq:mfe}, on every time interval, we have
\[
x(t)=y^{0,[n]}(t)+ O(\eps), \  v_\parallel(t)=v_\parallel^{0,[n]}(t)+ O(\eps), \quad n\eps \leq t \leq (n+1)\eps,
\]
and thus
\[
r(t)=r^{0,[n]} + O(\eps), \ z(t)=z^{0,[n]} + O(\eps), \ v_\parallel(t)=v_\parallel^{0,[n]} + O(\eps), \quad n\eps \leq t \leq (n+1)\eps.
\]

In view of the factor $\eps$ in front of the right hand side of the differential equations \eqref{eq:limit} and \eqref{eq:gc}, we have
\[
r^{0,[0]}(t)-\tilde{r}(t)=O(\eps), \quad z^{0,[0]}(t)-\tilde{z}(t)=O(\eps),  \quad v^{0,[0]}_\parallel(t)-\tilde{v}(t)=O(\eps), \quad  0\leq t \leq c/\eps.
\]
Since $y^{0,[n-1]}(n\eps)=y^{0,[n]}(n\eps)+O(\eps^N), \dot{y}^{0,[n-1]}(n\eps)=\dot{y}^{0,[n]}(n\eps)+O(\eps^{N-1})$, we have
\[
\begin{aligned}
r^{0,[n-1]}(n\eps)&=r^{0,[n]}(n\eps)+O(\eps^N)\\
z^{0,[n-1]}(n\eps)&=z^{0,[n]}(n\eps)+O(\eps^N)\\
v_\parallel^{0,[n-1]}(n\eps)&=v_\parallel^{0,[n]}(n\eps)+O(\eps^{N-1}).
\end{aligned}
\]
In view of the factor $\eps$ in front of the right hand side of the \eqref{eq:gc}, we have
\[
\begin{aligned}
r^{0,[n]}(t)-r^{0,[n-1]}(t)&=O(\eps^N), \\
z^{0,[n]}(t)-z^{0,[n-1]}(t)&=O(\eps^N),  \quad n\eps\leq t\leq c/\eps.\\
v_\parallel^{0,[n]}(t)-v_\parallel^{0,[n-1]}(t)&=O(\eps^{N-1}), 
 \end{aligned}
\]
With the above estimates, we obtain, for $n\eps \leq t\leq (n+1)\eps\leq c/\eps$
\[
\begin{aligned}
r(t)-\tilde{r}(t)&=r(t)-r^{0,[n]}+\sum_{j=1}^{n} \left(r^{0,[j]}(t)-r^{0,[j-1]}(t)\right)+ r^{0,[0]}(t)-\tilde{r}(t)\\
&=O(\eps)+O(n\eps^N)+O(\eps)=O(\eps), \\
z(t)-\tilde{z}(t)&=z(t)-z^{0,[n]}+\sum_{j=1}^{n} \left(z^{0,[j]}(t)-z^{0,[j-1]}(t)\right)+ z^{0,[0]}(t)-\tilde{z}(t)\\
&=O(\eps)+O(n\eps^N)+O(\eps)=O(\eps), \\
v_\parallel(t)-\tilde{v}(t)&=v_\parallel(t)-v_\parallel^{0,[n]}+\sum_{j=1}^{n} \left(v_\parallel^{0,[j]}(t)-v_\parallel^{0,[j-1]}(t)\right)+ v_\parallel^{0,[0]}(t)-\tilde{v}(t)\\
&=O(\eps)+O(n\eps^{N-1})+O(\eps)=O(\eps), \\
\end{aligned}
\]
which is the stated result of Theorem~\ref{thm:exact}.

\subsection{Proof of Theorem~\ref{thm:num}}
Similar to the proof of Theorem~\ref{thm:exact}, we structure the proof into three parts.

(a)
For a general strong magnetic field, the time interval of modulated Fourier expansion for numerical solution is validated over $O(h)$. Using the uniqueness of the modulated Fourier expansion, we can patch together many short-time expansions in the same way as it was done for the exact solution and obtain the expansion for longer time $O(1/\eps)$.

From Theorem 4.2 of \cite{lubich2022large}, it is known that the numerical solution $x^n$ given by the modified Boris algorithm \eqref{mboris}-\eqref{mod-init} with a step size $h$ satisfying
\[
c_*\eps \le h^2 \le C_*\eps
\]
can be written as
\begin{equation}\label{eq:mfe_num}
x^n=y^0(t_n) + (-1)^n y^1(t_n) +R_N(t_n), \qquad t_n=nh \le c/\eps,
\end{equation}
where $y^0=O(1)$, $y^1=O(h^2)$ are peicewise continuous with jumps of size $O(h^N)$ at integral multiples of $h$ and smooth elsewhere. They are unique up to $O(h^N)$ and $P_\perp(y^0)\dot{y}^0=O(h^2)$, $P_0(y^0)y^1=O(h^4)$.

Inserting \eqref{eq:mfe_num} into the numerical scheme \eqref{mboris} and separating the terms without $(-1)^n$ give the equation for  guiding centre $y^0(t)$ 
\begin{equation}\label{gc_Boris}
\begin{aligned}
\ddot{y}^0+h^2\ddddot{y}^0+O(h^4)=&\bigl(\dot{y}^0+ h^2\dddot{y}^0+ O(h^4)\bigr)\times B(y^0)+E(y^0)-\mu^0\,\nabla|B|(y^0)\\
&-\underbrace{\dot{y}^1_\perp\times B'(y^0)y^1_\perp}_{=:III}+O(h^4),
\end{aligned}
\end{equation}
where $III=O(h^2)$ in our stepsize regime $h^2\sim\eps$.
Taking the projection $P_{\pm 1}=P_{\pm 1}(y^0)$ on both sides gives
\[
\begin{aligned}
P_{\pm 1}\ddot{y}^0+O(h^2)=&\pm\iu|B(y^0)|P_{\pm 1}(\dot{y}^0+ h^2\dddot{y}^0+ O(h^4))\\
&+P_{\pm 1}\left(E(y^0)-\mu^0\,\nabla|B|(y^0)\right)+O(h^2),
\end{aligned}
\]
which means (recall that $h^2\sim\eps$)
\[
P_{\pm 1}\dot{y}^0=-h^2P_{\pm 1}\dddot{y}^0\mp\frac{\iu}{|B|}P_{\pm1}\ddot{y}^0\pm\frac{\iu}{|B|}P_{\pm1}\left(E(y^0)-\mu^0\,\nabla|B|(y^0)\right)+O(\eps h^2).
\]
Denoting $g_{\pm1}=P_{\pm1}\dot{y}^0$, we have $P_{\pm1}\ddot{y}^0=\dot{g}_{\pm1}-\dot{P}_{\pm1}\dot{y}^0$ and $P_{\pm1}\dddot{y}^0=\ddot{g}_{\pm1}-2\dot{P}_{\pm1}\ddot{y}^0-\ddot{P}_{\pm1}\dot{y}^0$. By differentiation and substitution the derivatives of $g_{\pm1}$ can be removed and we have 
\[
P_{\pm1}\dot{y}^0=h^2(2\dot{P}_{\pm1}\ddot{y}^0+\ddot{P}_{\pm1}\dot{y}^0)\pm\frac{\iu}{|B|}\dot{P}_{\pm1}\dot{y}^0\pm\frac{\iu}{|B|}P_{\pm1}\left(E(y^0)-\mu_0\nabla|B|(y^0)\right)+O(\eps h^2).
\]
Using that $\dot{P}_{1}+\dot{P}_{-1}+\dot{P}_{0}=0$, we obtain 
\begin{equation}\label{eq:y_Boris}
\begin{aligned}
\dot{y}^0
&=P_{\parallel}\dot{y}^0+P_{1}\dot{y}^0+P_{-1}\dot{y}^0\\
&=P_{\parallel}\dot{y}^0-h^2(2\dot{P}_\parallel\ddot{y}^0+\ddot{P}_\parallel \dot{y}^0)+\frac{1}{|B|}P_{\parallel}\dot{y}^0\times \frac{\d \e_\parallel}{\d t}+\frac{1}{|B|^2}\left(E-\mu^0\,\nabla |B|\right)\times B+O(\eps h^2),
\end{aligned}
\end{equation}
with $B, \nabla |B|, P_\parallel=\e_\parallel \e_\parallel^\top$ and $E$ evaluated at the guiding centre $y^0$.  Compared to the guiding centre equation \eqref{eq:y} of the exact solution, it is noticed that there are additional $O(h^2)$ terms in \eqref{eq:y_Boris}.

(b) Next we derive the guiding centre equation in toroidal geometry where $y^0(t)$ can be written as 
\[
y^0=r(y^0)\e_r(y^0)+z(y^0)\e_z=:r^0\e_r(y^0)+z^0\e_z.
\]

\noindent
--- Multiplying \eqref{eq:y_Boris} with $\e_r^\top=\e_r(y^0)^\top$ gives
\begin{equation}\label{eq:r_Boris}
\begin{aligned}
\e_r^\top \dot{y}^0 =&-2h^2 \e_r^\top \dot{P}_\parallel\ddot{y}^0-h^2\e_r^\top \ddot{P}_\parallel \dot{y}^0\\
&+\frac{\eps v^0_\parallel}{b(r,z)} \e_r^\top \left( \e_\parallel\times \frac{\d \e_\parallel}{\d t}\right)+\frac{\eps}{b(r,z)} \e_r^\top \left( (E-\mu^0\,\nabla b(r,z))\times \e_\parallel \right)+O(\eps h^2)
\end{aligned}
\end{equation}
with $v^0_\parallel:=\e_\parallel^\top  \dot{y}^0$.  Compared to \eqref{eq:r}, the only difference comes from the first two terms on  the right hand side which we calculate in the following.

Multiplying \eqref{gc_Boris} with $\e_\parallel^\top=\e_\parallel(y^0)^\top$ gives 
\[
\e_\parallel^\top \ddot{y}^0=O(h^2),
\]
 then the first term on the right hand side of \eqref{eq:r_Boris} is 
\[
-2h^2\e_r^\top \dot{P}_\parallel\ddot{y}^0=-2h^2\e_r^\top (\dot{\e}_\parallel \e_\parallel^\top \ddot{y}^0)=O(h^4).
\]
Using \eqref{eq:deriv}, the second term on the right hand side of \eqref{eq:r_Boris} can be expressed as
\begin{equation}\label{eq:formula}
\begin{aligned}
\e_r^\top\ddot{P}_\parallel \dot{y}^0&=\e_r^\top( \ddot{\e}_\parallel \e_\parallel^\top  \dot{y}^0)+2\e_r^\top( \dot{\e}_\parallel\dot{\e}_\parallel^\top \dot{y}^0)\\
&=-v^0_\parallel\frac{\d}{\d t}\left(\frac{v^0_\parallel}{r^0}\right)+2\left(\frac{v^0_\parallel}{r^0}\right)^2\frac{\d r^0}{\d t}\\
&=-\frac{v^0_\parallel}{r^0}\frac{\d v^0_\parallel}{\d t} + 3
\left(\frac{v^0_\parallel}{r^0}\right)^2\frac{\d r^0}{\d t}.
\end{aligned}
\end{equation}
Inserting 
\[
\frac{\d v^0_\parallel}{\d t}=\frac{\d}{\d t} (\e_\parallel^\top\dot{y}^0)= \dot{\e}_\parallel^\top \dot{y}^0+ \e_\parallel^\top\ddot{y}^0=-\frac{v^0_\parallel}{r^0}\frac{\d r^0}{\d t}+ O(h^2)
\]
into \eqref{eq:formula} gives
\[
\e_r^\top\ddot{P}_\parallel \dot{y}^0=4\left(\frac{v^0_\parallel}{r^0}\right)^2\frac{\d r^0}{\d t}+O(h^2).
\]
Then \eqref{eq:r_Boris} can be expressed as
\[
 \left(1+4h^2\left(\frac{v^0_\parallel}{r^0}\right)^2\right)\frac{\d r^0}{\d t}=\frac{\eps}{b} \e_r^\top \left((E-\mu^0\,\nabla b)\times \e_\parallel \right) +O(\eps h^2),
\]
which yields
\[
\frac{\d r^0}{\d t}=-\eps\frac{E_z}{b}+\eps\frac{\mu^0}{b}\partial_z b+O(\eps h^2),
\]
where the functions $E_z, b, \partial_z b$ are evaluated at $(r^0,z^0)$.
The initial value of $r^0$ can be expressed as
\[
r^0(0)=\langle y^0(0),\e_r(y(0))\rangle=\langle x^0,\e_r(x^0))\rangle + O(h^2)=r(x^0)+O(h^2).
\]

\noindent
--- Multiplying \eqref{eq:y_Boris} with $\e_z^\top$ gives
\begin{equation}\label{eq:z_Boris}
\begin{aligned}
\e_z^\top \dot{y}^0 =&-2h^2 \e_z^\top \dot{P}_\parallel\ddot{y}^0 -h^2 \e_z^\top \ddot{P}_\parallel \dot{y}^0\\
&+\frac{\eps v_\parallel}{b(r,z)} \e_z^\top\left( \e_\parallel\times \frac{\d \e_\parallel}{\d t}\right)+\frac{\eps}{b(r,z)} \e_z^\top \left((E-\mu^0\,\nabla b(r,z))\times \e_\parallel\right) +O(\eps h^2),
\end{aligned}
\end{equation}
where the first two terms on the right hand side vanish using \eqref{eq:deriv} and the orthorgonality of $\e_z,\e_r,\e_\parallel$. Similar to the continuous case, we obtain
\[
\frac{\d z^0}{\d t}=\eps\frac{\left(v^0_\parallel\right)^2}{r^0b}+\eps\frac{E_r}{b}-\eps\frac{\mu^0}{b}\partial_r b+O(\eps h^2),
\]
where the functions $E_r, b, \partial_r b$ are evaluated at $(r^0,z^0)$.
The initial value of $z^0$ is
\[
z^0(0)=\langle y^0(0),\e_z)\rangle=\langle x^0,\e_z)\rangle + O(h^2)=z(x^0)+O(h^2).
\]

\noindent
--- Finally, we need to derive the differential equation for $v^0_\parallel$ which can be directly computed as the continuous case
\begin{equation}\label{eq:parallel_Boris}
\frac{\d v^0_\parallel}{\d t}=\frac{\d}{\d t}\left( \e_\parallel^\top \dot{y}^0\right)= \dot{\e}_\parallel^\top \dot{y}^0+ \e_\parallel^\top\ddot{y}^0.
\end{equation}
The first term on the right hand side is the same as \eqref{eq:form}.
Multiplying \eqref{gc_Boris} with $\e_\parallel^\top=\e_\parallel(y^0)^\top$ yields
\begin{equation}\label{eq:ddoty_Boris}
\begin{aligned}
 \e_\parallel^\top \ddot{y}^0&=-h^2 \e_\parallel^\top \ddddot{y}^0+ O(h^4)\\
&=-h^2 \frac{\d^2}{\d t^2}(\e_\parallel^\top \ddot{y}^0)+2h^2 (\dot{\e}_\parallel^\top \dddot{y}^0)+h^2 \ddot{\e}_\parallel^\top\ddot{y}^0+ O(h^4),
\end{aligned}
\end{equation}
where we use $ \e_\parallel^\top(E-\mu^0\,\nabla|B|) =0$ and $\e_\parallel^\top III =0$. 
Since the derivatives of $r^0$ are $O(\eps)$ and using \eqref{eq:deriv}, we have
\[
\begin{aligned}
 \dot{\e}_\parallel^\top \dddot{y}^0&=-\frac{v^0_\parallel}{r^0}  \e_r^\top \dddot{y}^0\\
&=-\frac{v^0_\parallel}{r^0}\left( \dddot{r^0} - 3 \dot{\e}_r^\top \ddot{y}^0 - 3 \ddot{\e}_r^\top\dot{y}^0- \dddot{\e}_r^\top y^0 \right)\\
&=-\frac{v^0_\parallel}{r^0}\left(- 3\frac{v^0_\parallel}{r^0} \e_\parallel^\top \ddot{y}^0 -3v^0_\parallel\frac{\d}{\d t}\left(\frac{v^0_\parallel}{r^0}\right) +  v^0_\parallel\frac{\d}{\d t}\left(\frac{v^0_\parallel}{r^0}\right) + r^0\frac{\d}{\d t}\left(\left(\frac{v^0_\parallel}{r^0}\right)^2\right) \right)+O(\eps)\\
&=3\left(\frac{v^0_\parallel}{r^0}\right)^2 \e_\parallel^\top \ddot{y}^0+O(\eps)
\end{aligned}
\]
and
\[
\begin{aligned}
 \ddot{\e}_\parallel^\top \ddot{y}^0&=-\e_r^\top \ddot{y}^0\frac{\d}{\d t}\left(\frac{v^0_\parallel}{r^0}\right)-\left(\frac{v^0_\parallel}{r^0}\right)^2 \e_\parallel^\top \ddot{y}^0\\
 &=-\left(\frac{\d}{\d t}(\e_r^\top \dot{y}^0)-\dot{\e}_r^\top\dot{y}^0\right)\frac{\d}{\d t}\left(\frac{v^0_\parallel}{r^0}\right)-\left(\frac{v^0_\parallel}{r^0}\right)^2 \e_\parallel^\top\ddot{y}^0\\
 &=\frac{\left(v^0_\parallel\right)^2}{r^0}\frac{\d}{\d t}\left(\frac{v^0_\parallel}{r^0}\right)-\left(\frac{v^0_\parallel}{r^0}\right)^2 \e_\parallel^\top\ddot{y}^0+O(\eps),
 \end{aligned}
\]
then \eqref{eq:ddoty_Boris} can be written as
\[
 (1-5h^2(v^0_\parallel/r^0)^2) \e_\parallel^\top \ddot{y}^0=h^2\left(\frac{v^0_\parallel}{r^0}\right)^2\frac{\d v^0_\parallel}{\d t}+ O(h^4).
\]
This gives
\[
 \e_\parallel^\top \ddot{y}^0=h^2\left(\frac{v^0_\parallel}{r^0}\right)^2\frac{\d v^0_\parallel}{\d t}+ O(h^4).
\]
\eqref{eq:parallel_Boris} now can be written as
\[
(1-h^2(v^0_\parallel/r^0)^2)\frac{\d v^0_\parallel}{\d t}=-\frac{v^0_\parallel}{r^0}\frac{\d r^0}{\d t}+ O(h^4),
\]
which gives
\[
\frac{\d v^0_\parallel}{\d t}=-\frac{v^0_\parallel}{r^0}\frac{\d r^0}{\d t}+ O(h^4).
\]
The initial value of $v^0_\parallel$ is
\[
v^0_\parallel(0)=\langle \dot{y}^0(0), \e_\parallel(y^0(0))\rangle=\langle v^0,\e_\parallel(x^0) \rangle + O(h^2).
\]

(c)
Denoting by $y^{0,[n]}, r^{0,[n]}, z^{0,[n]}, v^{0,[n]}_\parallel$ the functions  $y^{0}, r^0, z^0, v^0_\parallel$  on time interval $nh\leq t\leq (n+1)h$, from (b), it is known that  these coefficients satisfy the following equations
\[
\begin{aligned}
\frac{\d r^{0,[n]}}{\d t}&=-\eps\frac{E_z}{b}+\eps\frac{\mu^0}{b}\partial_z b+O(\eps h^2), \quad r^{0,[n]}(nh)=r(x(nh))+O(h^2)\\
\frac{\d z^{0,[n]}}{\d t}&=\eps\frac{(v^{0,[n]}_\parallel)^2}{b r^{0,[n]}}+\eps\frac{E_r}{b}-\eps\frac{\mu^0}{b}\partial_r b+O(\eps h^2), \quad 
z^{0,[n]}(nh)=z(x(nh))+O(h^2)\\
\frac{\d v^{0,[n]}_\parallel}{\d t}&=\eps\frac{v^{0,[n]}_\parallel}{r^{0,[n]}}\left(\frac{E_z}{b}-\frac{\mu^0}{b}\partial_z b\right)+O(\eps h^2), \quad v^{0,[n]}_\parallel(nh)=\langle \dot{x}(nh), \e_\parallel(x(nh)) \rangle + O(h^2),
\end{aligned}
\]
with $E_r, E_z, b, \partial_r b, \partial_z b$ evaluated at $(r^{0,[n]},z^{0,[n]})$.

By patching together the errors together as what was done for the continuous case, we prove that
\[
r(x^n)-\tilde{r}(t_n)=O(h^2), \quad z(x^n)-\tilde{z}(t_n)=O(h^2),  \quad v^n_\parallel-\tilde{v}(t_n)=O(h^2), \quad  0\leq t \leq c/\eps.
\]


\section*{Acknowledgement} 
The author thanks Professor Christian Lubich for many useful discussions and comments.
This work was supported by the Sino-German (CSC-DAAD) Postdoc Scholarship, Program No. 57575640.

\bibliographystyle{abbrv}
\bibliography{toroidal}

\end{document}